\documentclass{article}

\newcommand{\be}{\begin{equation}}
\newcommand{\ee}{\end{equation}}
\newcommand{\bea}{\begin{eqnarray}}
\newcommand{\eea}{\end{eqnarray}}
\newcommand{\barray}{\begin{array}}
\newcommand{\earray}{\end{array}}
\newcommand{\pa}{\partial}
\newcommand{\nn}{\nonumber}
\newcommand{\bitem}{\begin{itemize}}
\newcommand{\eitem}{\end{itemize}}
\newtheorem{teo}{Theorem}[section]
\newcommand{\bt}{\begin{teo}}
\newcommand{\et}{\end{teo}}
\newtheorem{Def}{Definition}[section]
\newcommand{\bd}{\begin{Def}}
\newcommand{\ed}{\end{Def}}
\newtheorem{lem}{Lemma}[section]
\newcommand{\bl}{\begin{lem}}
\newcommand{\el}{\end{lem}}
\newtheorem{prop}{Proposition}[section]
\newcommand{\bp}{\begin{prop}}
\newcommand{\ep}{\end{prop}}
\newtheorem{cor}{Corollary}[section]
\newcommand{\bc}{\begin{cor}}
\newcommand{\ec}{\end{cor}}
\newtheorem{ex}{Example}[section]
\newcommand{\bex}{\begin{ex}}
\newcommand{\eex}{\end{ex}}
\newtheorem{rem}{Remark}[section]
\newcommand{\br}{\begin{rem}}
\newcommand{\er}{\end{rem}}

\catcode `\@=11
\@addtoreset{equation}{section}

\begin{document}

\begin{center}
{\Large \textbf{Compatible Dubrovin--Novikov \\ Hamiltonian
operators, Lie derivative and \\ integrable
systems of hydrodynamic type\footnote{This work was supported by
the Alexander von Humboldt Foundation (Germany),
the Russian Foundation for Basic Research
(grant No. 99--01--00010) and the INTAS
(grant No. 99--1782). The work will be published in
Proceedings of the International Conference ``Nonlinear
Evolution Equations and Dynamical Systems,'' Cambridge (England),
July 24--30, 2001.}}}
\end{center}

\bigskip
\bigskip

\centerline{\large {O. I. Mokhov}}
\bigskip
\bigskip

\section{Introduction}

In the present paper, we prove that a local Hamiltonian
operator of hydrodynamic type $K^{ij}_1$ (Dubrovin--Novikov
Hamiltonian operator \cite{[1]}) is compatible with
a nondegenerate local Hamiltonian operator of
hydrodynamic type $K^{ij}_2$ if and only if the operator
$K^{ij}_1$ is locally the Lie derivative of the
operator $K^{ij}_2$ along a vector field
in the corresponding domain of local coordinates.
This result gives, first of all, a convenient
general invariant criterion of the compatibility
for the Dubrovin--Novikov Hamiltonian operators and,
in addition, this gives a natural invariant
definition of the class of special flat manifolds
corresponding to all the class of compatible
Dubrovin--Novikov Hamiltonian operators
(the Frobenius--Dubrovin manifolds naturally belong
to this class of flat manifolds).
There is an integrable bi-Hamiltonian hierarchy
corresponding to every flat manifold of this class.
The integrable systems are also studied in
the present paper. This class of integrable systems
is explicitly given by solutions of the nonlinear
system of equations, which is integrated by
the method of inverse scattering problem.
The corresponding results on compatible
nonlocal Hamiltonian operators of
hydrodynamic type and other Hamiltonian and
symplectic differential-geometric type operators
related by the Lie derivative and the results on
the corresponding to them integrable bi-Hamiltonian systems
will be published in other our works.

Recall that an operator $K^{ij}[u(x)]$
is called {\it Hamiltonian} if it defines
a Poisson bracket
(skew-symmetric and satisfying the Jacobi identity)
\be
\{ I, J \} = \int {\delta I \over \delta u^i (x)}
K^{ij} [u(x)] {\delta J \over \delta u^j (x)} dx \label{po}
\ee
for arbitrary functionals
$I[u(x)]$ and $J[u(x)]$ on the space of functions
(fields) $u(x) = \{ u^i(x),\ 1\leq i\leq N \},$ where
$u^1,...,u^N$ are local coordinates on a certain
given smooth $N$-dimensional manifold $M$.
It is obvious that any Hamiltonian operator
always behaves as a contravariant two-valent operator
tensor field with respect to arbitrary local changes of
coordinates
$u^i = u^i (v^1,...,v^N), \ 1 \leq i \leq N,$
on the manifold $M$:
\be
\widetilde{K}^{sr}[v(x)] = {\pa v^s \over \pa u^i}
K^{ij} [u(v(x))] \circ {\pa v^r \over \pa u^j}, \label{ten}
\ee
where the symbol $\circ $ means the
operator multiplication
(it is important to note
that in contrast to the usual tensor fields
the components of operators in (\ref{ten}) multiply only
in the indicated order).
Accordingly, the classical tensor constructions,
in particular, the Lie derivative along a vector field,
are applied to Hamiltonian operators
(see, for example, the monograph
\cite{[2]} and also section \ref{sect2} below).

Hamiltonian operators are called {\it compatible}
if any their linear combination is also a Hamiltonian operator
(Magri, \cite{[3]}; see also \cite{[2]}).
In the theory of Hamiltonian operators, there is
the following general fact (see, for example,
\cite{[2]}), which is important for applications:
if the second cohomology group of the corresponding
complex of formal variational calculus is trivial,
then from the compatibility of two Hamiltonian
operators $K^{ij}_1$ and $K^{ij}_2$, where
$K^{ij}_2$ is an invertible operator,
it follows that there exists
a formal vector field $X[u(x)]$ depending, generally
speaking, in an arbitrary and nonlocal way,
on the fields $u(x)$ and their derivatives and such that
$K_1 = L_X K_2,$ where $L_X K$
is the Lie derivative of a Hamiltonian operator $K$
along a formal vector field $X$.
Moreover, if an operator $K$ is Hamiltonian and,
in addition, the operator $L_X K$, where $X$
is a certain arbitrary formal vector field,
is also Hamiltonian, then the Hamiltonian operators
$K$ and $L_X K$ are always compatible.
If an operator $K$ is Hamiltonian and, besides, $L^2_X K = 0$,
where $X$ is a certain vector field,
then the operator $L_X K$ is always Hamiltonian,
and consequently the Hamiltonian operators
$K$ and $L_X K$ are compatible in this case.
This beautiful construction plays the very important
role in applications and the
corresponding class of compatible Hamiltonian
operators (the pairs of operators of the form
$K^{ij},$ $(L_X K)^{ij},$ such that  $(L^2_X K)^{ij} = 0$)
deserves a separate study.

An important special partial case
of this construction arose afterwards in the Dubrovin
theory
of Frobenius manifolds (the quasihomogeneous
compatible nondegenerate Dubrovin--Novikov Hamiltonian
operators) \cite{[4]}--\cite{[6]}.
The study of general compatible Dubrovin--Novikov
Hamiltonian operators of this class, that is,
of the form $B^{ij}$ and $(L_X B)^{ij}$ and such that
$L_X^2 B = 0,$ was started by the present
author and Fordy in \cite{[7]}, where the first
classification results were obtained.

\section{Compatible Dubrovin--Novikov \\ Hamiltonian
operators}  \label{sect1}

Recall that a Hamiltonian operator given by
an arbitrary matrix homogeneous first-order
ordinary differential operator, that is,
a Hamiltonian operator of the form
\be
P^{ij} [u (x)] =
g^{ij} (u(x))\, {d \over dx} + b^{ij}_k (u(x)) \, u^k_x,
\label{hydro}
\ee
is called {\it a local Hamiltonian operator of hydrodynamic type}
or {\it Dubrovin--Novikov Hamiltonian operator} [1].
This definition does not depend on the choice
of local coordinates $u^1,...,u^N$ on the manifold
$M$, since it follows from (\ref{ten}) that
the form of operator (\ref{hydro}) is invariant
with respect to local changes of coordinates on
$M$. Operator (\ref{hydro}) is called {\it nondegenerate}
if $\det (g^{ij} (u)) \not\equiv 0.$
If $\det (g^{ij} (u)) \not\equiv 0,$ then
operator (\ref{hydro}) is Hamiltonian if and only if
1) $g^{ij} (u)$ is an arbitrary contravariant flat
pseudo-Riemannian metric (a metric of zero
Riemannian curvature),
2) $b^{ij}_k (u) = - g^{is} (u) \Gamma ^j_{sk} (u),$ where
$\Gamma^j_{sk} (u)$ is the Levi-Civita
connection generated by the metric $g^{ij} (u)$
(the Dubrovin--Novikov theorem \cite{[1]}).
In particular, it follows from the Dubrovin--Novikov
theorem that for any nondegenerate local Hamiltonian
operator of hydroodynamic type there always exist local
coordinates
$v^1,...,v^N$ (flat coordinates of the metric $g^{ij}(u)$)
in which all the coefficients of the operator are constant:
\be
\widetilde g^{ij} (v)
= \eta^{ij} = {\rm \ const}, \ \
\widetilde \Gamma^i_{jk} (v) = 0, \ \
\widetilde b^{ij}_k (v) =0,
\ee
that is the corresponding Poisson bracket
has the form
\be
\{ I,J \} = \int {\delta I \over \delta v^i(x)}
 \eta^{ij} {d \over dx}
{\delta J \over \delta v^j(x)} dx,
\label{const}
\ee
where $(\eta^{ij})$ is a nondegenerate
symmetric constant matrix:
\be
\eta^{ij} = \eta^{ji}, \ \ \eta^{ij} = {\rm const},
\  \ \det \, (\eta^{ij}) \neq 0.\nn
\ee

Moreover, it immediately follows from the Dubrovin--Novikov
theorem that any two nondegenerate Dubrovin--Novikov
Hamiltonian operators
$P^{ij}_1 [u (x)]$ and $P^{ij}_2 [u(x)]$ generated
by flat contravariant metrics $g^{ij}_1 (u)$ and
$g^{ij}_2 (u)$ respectively are compatible
if and only if
1) any linear combination of these flat metrics
\be
g^{ij} (u) = \lambda_1 g_1^{ij} (u) + \lambda_2 g_2^{ij} (u),
\label{comb0}
\ee
where $\lambda_1$ and $\lambda_2$ are arbitrary
constants for which
$\det ( g^{ij} (u) ) \not\equiv 0$,
is also a flat metric,
2) the coefficients of the corresponding
Levi-Civita connections are related by the same
linear formula:
\be
\Gamma^{ij}_k (u) = \lambda_1 \Gamma^{ij}_{1, k} (u) +
\lambda_2 \Gamma^{ij}_{2, k} (u). \label{sv0}
\ee
The derived purely differential-geometric
conditions (\ref{comb0}) and
(\ref{sv0}) on flat metrics
$g_1^{ij} (u)$ and $g_2^{ij} (u)$ define
{\it a flat pencil of metrics} \cite{[4]}.
In this case, we shall also say that the flat metrics
$g_1^{ij} (u)$ and $g_2^{ij} (u)$
are {\it compatible} (see \cite{[8]}).

So the problem of description for
compatible nondegenerate local Hamiltonian operators of
hydrodynamic type
is the purely differential-geometric problem of
description of general flat pencils of metrics
(see \cite{[4]}).
In \cite{[4]}, \cite{[5]} Dubrovin considered all the
tensor relations for the general flat pencils of metrics.

First of all, let us introduce the necessary notation.
Let $\nabla_1$ and $\nabla_2$ be the operators of
covariant differentiation given by the Levi-Civita
connections $\Gamma^{ij}_{1,k} (u)$ and $\Gamma^{ij}_{2,k} (u)$
generated by the metrics $g^{ij}_1 (u)$ and $g^{ij}_2 (u)$
respectively. The indices of the covariant
differentials are raised and lowered by the
corresponding metrics:
$\nabla^i_1= g^{is}_1 (u) \nabla_{1,s}$,
$\nabla^i_2=g^{is}_2 (u) \nabla_{2,s}$.
Consider the tensor
\be
\Delta ^{ijk} (u) = g^{is}_1 (u) g^{jp}_2 (u)
\left (\Gamma^k_{2, ps} (u)
- \Gamma^k_{1, ps} (u) \right ),  \label{(2.3)}
\ee
introduced by Dubrovin in \cite{[4]}, \cite{[5]}.

\bt [Dubrovin \cite{[4]}, \cite{[5]}] \label{dub1}
If metrics $g^{ij}_1 (u)$ and $g^{ij}_2 (u)$ form
a flat pencil, then there exists a vector field
$f^i (u)$ such that the tensor
$\Delta ^{ijk} (u)$ and the metric $g^{ij}_1 (u)$
have the form
\be
\Delta ^{ijk} (u) = \nabla^i_2 \nabla^j_2 f^k (u),
\label{(2.4)}
\ee
\be
g^{ij}_1 (u) = \nabla^i_2 f^j (u) + \nabla^j_2 f^i (u) +
c g^{ij}_2 (u), \label{(2.5)}
\ee
where $c$ is a certain constant, and the vector field
$f^i (u)$ satisfies the equations
\be
\Delta^{ij}_s (u) \Delta^{sk}_l (u) =
\Delta^{ik}_s (u) \Delta^{sj}_l (u), \label{(2.6)}
\ee
where
\be
\Delta^{ij}_k (u) =g_{2,ks} (u) \Delta ^{sij} (u)
= \nabla_{2,k} \nabla^i_2 f^j (u), \label{(2.7)}
\ee
and
\be
(g^{is}_1 (u) g^{jp}_2 (u) - g^{is}_2 (u) g^{jp}_1 (u))
\nabla_{2,s} \nabla_{2,p} f^k (u) =0. \label{(2.8)}
\ee
Conversely, for the flat metric $g^{ij}_2 (u)$
and the vector field $f^i (u)$ that is a solution
of the system of equations (\ref{(2.6)}) and (\ref{(2.8)}),
the metrics $g^{ij}_2 (u)$ and (\ref{(2.5)})
form a flat pencil.
\et

In \cite{[6]} Dubrovin proved that
the theory of Frobenius manifolds
constructed him in \cite{[4]} (the Frobenius manifolds
correspond to two-dimensional topological
field theories) is equivalent to the theory
of quasihomogeneous compatible nondegenerate
Dubrovin--Novikov Hamiltonian operators
or, in other words, quasihomogeneous
flat pencils of metrics.

A flat pencil of metrics
generated by flat metrics $g^{ij}_1 (u)$ and $g^{ij}_2 (u)$
is called {\it quasihomogeneous of degree $d$} if
there exists a function $\tau (u)$ such that
for the vector fields
\be
E = \nabla_1 \tau, \ \  E^i = g^{is}_1 (u) {\partial \tau \over
\partial u^s}, \ \ \ \ e = \nabla_2 \tau, \ \ e^i =
g^{is}_2 (u) {\partial \tau \over
\partial u^s},
\ee
the following conditions are satisfied:

1) $[e, E] = e,$

2) $L_E g_1 = (d-1)g_1,$

3) $L_e g_1 = g_2,$

4) $L_e g_2 = 0,$

\noindent
where $L_X g$ is the Lie derivative of a metric
$g$ along a vector field
$X$ (Dubrovin, \cite{[6]}).
Note that in the quasihomogeneous case we always have
$L_e Q_1 = Q_2,$
$L_e Q_2 = 0,$ $L^2_e Q_1 =0,$ for the corresponding
quasihomogeneous compatible nondegenerate Dubrovin--Novikov
Hamiltonian operators
$Q^{ij}_1 [u(x)]$ and $Q^{ij}_2 [u(x)]$ (more in detail
about this important class of quasihomogeneous compatible
nondegenerate Dubrovin--Novikov Hamiltonian operators
see \cite{[7]}).

In the present author's work \cite{[9]},
a necessary for further, explicit and
simple {\em ctiterion} of compatibility
for two local Poisson brackets of hydrodynamic type
was stated in a suitable for us way, that is,
it is shown what an explicit form is sufficient and
necessary for two local Hamiltonian operators of
hydrodynamic type to be compatible (see also theorem \ref{dub1}).

\bl [\cite{[9]}]  \label{lem1}
{\bf (explicit criterion of compatibility
for local Poisson brackets of hydrodynamic type)}
Any local Poisson bracket of hydrodynamic type
$\{ I, J \}_2$  is compatible with the constant
nondegenerate Poisson bracket (\ref{const})
if and only if it has the form
\be
\{ I, J \}_2 =
\int {\delta I \over \delta v^i (x)} \biggl (
\biggl ( \eta^{is} {\partial h^j \over \partial v^s}
+ \eta^{js} {\partial h^i \over \partial v^s} \biggr )
{d \over dx} + \eta^{is}
{\partial^2 h^j \over \partial v^s \partial v^k}
v^k_x \biggr ) {\delta J \over \delta v^j (x)} dx, \label{(2.9)}
\ee
where $h^i (v), \ 1 \leq i \leq N,$ are smooth functions
defined in a certain domain of local coordinates.
\el

We do not require in lemma \ref{lem1}
that the Poisson bracket of hydrodynamic type
$\{ I, J \}_2$ is nondegenerate.
Besides, it is important to note that this statement is local.

\section{Lie derivative and Dubrovin--Novikov \\
Hamiltonian operators} \label{sect2}

Let $K^{ij}[u(x)]$ be an arbitrary Hamiltonian operator,
$\xi (u) = \{\xi ^i (u),\ 1 \leq i \leq N \}$
is an arbitrary smooth vector field on the manifold $M$.
{\it The Lie derivative of the Hamiltonian
operator $K^{ij}[u(x)]$}
(just as any contravariant two-valent operator tensor
field of type (\ref{ten}))
{\it along the vector field $\xi(u)$} is the operator
\be
({\mathcal{L}}_{\xi}K)^{ij} [u(x)] =
\left. \left (
{d \over dt} \left [ \left ( \delta^i_s - t {\pa \xi^i \over \pa u^s}
\right ) K^{sr}[u(x) + t \xi (u(x))] \circ \left ( \delta^j_r - t
{\pa \xi^j \over \pa u^r} \right ) \right ] \right ) \right |_{t=0},
\label{flie}
\ee
which also always behaves as contravariant two-valent
operator tensor field of type (\ref{ten})
with respect to arbitrary local changes of
coordinates on the manifold $M$.

For the Lie derivative of a local Hamiltonian operator
of hydrodynamic type $P^{ij}[u(x)]$ along the vector
field $\xi (u)$ we get
\bea
&&
({\mathcal{L}}_{\xi} P)^{ij} [u(x)]=
\left ( \xi^s  {\pa g^{ij} \over \pa u^s} - g^{sj}
{\pa \xi^i \over \pa u^s} - g^{is} {\pa \xi^j \over \pa u^s} \right )
{d \over dx} \nn\\
&&
+ \left ( \xi^s {\pa b^{ij}_k \over \pa u^s}
-b^{sj}_k {\pa \xi^i \over \pa u^s} - b^{is}_k {\pa \xi ^j \over \pa u^s}
+ b^{ij}_s {\pa \xi ^s \over \pa u^k} -
g^{is} {\pa^2 \xi^j \over \pa u^s \pa u^k } \right ) u^k_x.
\label{liehyd}
\eea

\bt \label{prlie}
{\it Any Dubrovin--Novikov Hamiltonian
operator $P^{ij}_1$ is compatible with
a nondegenerate Dubrovin--Novikov Hamiltonian
operator $P^{ij}_2$ if and only if there locally
exists a vector field $\xi (u)$ such that
\be
P^{ij}_1 = ({\mathcal{L}}_{\xi} P_2)^{ij}. \label{li}
\ee
}
\et

In a different form, which is not connected with
the Lie derivative, the general compatibility conditions
for the Dubrovin--Novikov Hamiltonian operators were
shown in \cite{[4]}, \cite{[8]}--\cite{[11]}
 (see also section \ref{sect1} above).

{\large Question.} If there always globally exists
such a smooth vector field on every flat manifold $M$,
on the loop space of which are globally defined
compatible Dubrovin--Novikov
Hamiltonian operators $P_1^{ij}$ and $P_2^{ij}$?
If not always, then it is necessary to investigate
the corresponding differential-geometric and topology
obstructions and also to single out and study all
the flat manifolds on which there globally exists
a vector field such that the Lie derivative of the
Dubrovin--Novikov Hamiltonian operator generated by
the flat metric of the manifold along the vector field
is also a Hamiltonian operator.
In particular, it is interesting to study
all such vector fields in ${\bf R}^N$ or ${\bf R}_{k, N-k}^N$,
and also in domains of these spaces. Locally, any
such vector field is a solution of the nonlinear
system of equations integrable by the method of
inverse scattering problem
(see \cite{[10]}--\cite{[11]}) and generates
integrable bi-Hamiltonian systems of
hydrodynamic type (we shall consider them in the
next section).

Here let us prove theorem \ref{prlie}.
Let Dubrovin--Novikov Hamiltonian operators
$P^{ij}_1 [u(x)]$ and $P^{ij}_2 [u(x)]$ be
compatible and the operator $P^{ij}_2 [u(x)]$ be
nondegenerate. Consider the local coordinates
$v=(v^1,...,v^N)$ in which the nondegenerate
Dubrovin--Novikov Hamiltonian operator
$P^{ij}_2 [u(x)]$ is reduced to the constant form
(\ref{const}):
\be
P^{ij}_2 [v(x)] = \eta^{ij} {d \over dx}, \label{opconst}
\ee
where ($\eta^{ij}$) is an arbitrary nondegenerate
constant symmetric matrix:
$\det (\eta^{ij}) \neq 0,$ $ \eta^{ij} = {\rm const},$
$\eta^{ij} = \eta^{ji}$ (there always exist
such local coordinates by the Dubrovin--Novikov theorem).
In these coordinates, according to formula
(\ref{liehyd}), for the Lie derivative of the operator
$P^{ij}_2 [v(x)]$ along an arbitrary vector field
$\xi (v) = (\xi^1(v),...,\xi^N(v))$ we get
\be
({\mathcal{L}}_{\xi} P_2)^{ij} [v(x)]=
\left (  - \eta^{sj}
{\pa \xi^i \over \pa v^s} - \eta^{is} {\pa \xi^j \over \pa v^s} \right )
{d \over dx} -
\eta^{is} {\pa^2 \xi^j \over \pa v^s \pa v^k } v^k_x.
\label{liehydconst}
\ee
According to lemma \ref{lem1} any Dubrovin--Novikov
Hamiltonian operator $P^{ij}_1 [v(x)]$
compatible with the Hamiltonian operator
(\ref{opconst}) must have namely such the form
(\ref{liehydconst}) in the local coordinates
$v=(v^1,...,v^N)$ (in formula
(\ref{(2.9)}) $h^i (v) = - \xi^i (v),$ $1 \leq i \leq N$).
Thus there exists a vector field  ${\xi}^i (v),$
$1 \leq i \leq N,$ such that
\be
P^{ij}_1 [v(x)] = ({\mathcal{L}}_{\xi} P_2)^{ij} [v(x)]. \label{li2}
\ee
Then by virtue of tensor invariance of the Lie derivative
formula (\ref{li}) is valid in arbitrary
local coordinates.

Conversely, let two Dubrovin--Novikov Hamiltonian operators
$P^{ij}_1 [u(x)]$ and $P^{ij}_2 [u(x)]$ are related by
formula
(\ref{li}) and the operator $P^{ij}_2 [u(x)]$
is nondegenerate. Then, in the local coordinates
$v=(v^1,...,v^N)$ in which the nondegenerate
Dubrovin--Novikov Hamiltonian operator $P^{ij}_2 [u(x)]$
is reduced to the constant form
(\ref{opconst}), the Hamiltonian operator
$P^{ij}_1 [v(x)]$ has the form (\ref{liehydconst}).
According to lemma \ref{lem1} a pair of Hamiltonian
operators
$P^{ij}_1 [v(x)]$ and $P^{ij}_2 [v(x)]$
of such the form are necessarily compatible. Thus
theorem \ref{prlie} is proved.

Now the special class of flat manifolds which corresponds
to the class of all compatible Dubrovin--Novikov
Hamiltonian operators and generalizes the class of
Frobenius--Dubrovin manifolds is naturally singled out.
Theorem \ref{prlie} gives the following invariant
definition of these manifolds.
We consider the manifolds $(M, g, \xi)$, where
$M$ is a flat manifold with a flat metric $g$
equipped also with a vector field $\xi$ such that,
for the nondegenerate Dubrovin--Novikov Hamiltonian
operator $P^{ij} [u(x)]$ generated by the flat metric
$g$, the operator
$({\mathcal{L}}_{\xi} P)^{ij} [u(x)]$ is also Hamiltonian
(that is, the operator defines
a Poisson bracket on the corresponding
loop space of the manifold $M$).
Generally speaking, for the description
of all compatible Dubrovin--Novikov Hamiltonian
operators it is necessary to consider the following more weak condition:
an existence of such vector field locally, in a neighbourhood
of every point of the manifold. Locally,
such vector fields are described by the nonlinear system
integrable by the method of inverse scattering problem
(see \cite{[10]}--\cite{[11]}).

\section{Class of integrable bi-Hamiltonian
systems \\ of hydrodynamic type} \label{sect3}

An arbitrary pair of compatible Dubrovin--Novikov
Hamiltonian operators $P^{ij}_1$ and $P^{ij}_2$, one
of which (let us assume $P^{ij}_2$) is nondegenerate,
can be reduced to the following special form
by a local change of coordinates:

\be
P^{ij}_2 [v(x)] = \eta^{ij} {d \over dx}, \label{op1}
\ee
\be
P^{ij}_1 [v(x)] =
\biggl ( \eta^{is} {\partial h^j \over \partial v^s}
+ \eta^{js} {\partial h^i \over \partial v^s} \biggr )
{d \over dx} + \eta^{is}
{\partial^2 h^j \over \partial v^s \partial v^k}
v^k_x,   \label{op2}
\ee
where ($\eta^{ij}$) is an arbitrary nondegenerate
constant symmetric matrix:
$\det (\eta^{ij}) \neq 0,$ $ \eta^{ij} = {\rm const},$
$\eta^{ij} = \eta^{ji};$  $h^i (v),$  $1 \leq i \leq N,$
are smooth functions given in a certain
domain of local coordinates such that operator (\ref{op2})
is Hamiltonian
(lemma \ref{lem1}, see also theorem \ref{dub1} above).
An operator of form
(\ref{op2}) is Hamiltonian if and only if
\be
\eta^{sr} {\partial^2 h^j	 \over \partial v^s \partial v^i}
 {\partial^2 h^k \over \partial v^l \partial v^r} =
\eta^{sr} {\partial^2 h^k	 \over \partial v^s \partial v^i}
 {\partial^2 h^j \over \partial v^l \partial v^r},
\label{(2.17)}
\ee
\be
\biggl ( \eta^{ip} {\partial h^s \over \partial v^p} +
\eta^{sp} {\partial h^i \over \partial v^p} \biggr )\eta^{jl}
 {\partial^2 h^r \over \partial v^l \partial v^s} =
\biggl (\eta^{jp} {\partial h^s \over \partial v^p} +
\eta^{sp} {\partial h^j \over \partial v^p} \biggr ) \eta^{il}
 {\partial^2 h^r \over \partial v^l \partial v^s}.
\label{(2.18)}
\ee
(\cite{[12]}, see also \cite{[9]} and theorem \ref{dub1}).

The system of nonlinear equations (\ref{(2.17)}), (\ref{(2.18)}),
as was conjectured in the present author's work \cite{[12]}
(there was stated the corresponding conjecture in \cite{[12]}),
is integrable by the method of inverse scattering problem.
The procedure of integrating this system was presented
in the author's work
\cite{[10]}, \cite{[10a]}.
In the work \cite{[11]} the Lax pair for this
system was demonstrated.
Note that the associativity equations of two-dimensional
topological field theory (see \cite{[4]}) are
a natural reduction of equations (\ref{(2.17)})
for ``potential'' vector fields $h^i (v),$ $1 \leq i \leq N,$
of the special form
\be
h^i (v) = \eta^{ij} {\partial \Phi \over \partial v^j},
\ee
where $\Phi (v)$ is a certain smooth function (`` the potential'')
(see also \cite{[9]}, \cite{[12]}--\cite{[14]}).

Consider the recursion operator
generated by the
``canonical'' compatible
Dub\-ro\-vin--No\-vi\-kov Hamiltonian
operators (\ref{op1}),
(\ref{op2}):

\be
R^i_l =
\left [ P_1 [v (x)] \left ( P_2 [v(x)] \right )^{-1} \right ]^i_l
=   \left ( \biggl ( \eta^{is} {\partial h^j \over \partial v^s}
+ \eta^{js} {\partial h^i \over \partial v^s} \biggr )
{d \over dx} + \eta^{is}
{\partial^2 h^j \over \partial v^s \partial v^k}
v^k_x \right ) \eta_{jl} \left ( {d \over dx} \right )^{-1},
\label{recur}
\ee
where ($\eta_{ij}$) is the matrix which is inverse
to the matrix
($\eta^{ij}$):\  $\eta^{is} \eta_{sj} = \delta^i_j$
(see \cite{[2]}, \cite{[15]}--\cite{[19]}
about recursion operators generated by pairs of
compatible Hamiltonian operators).

Let us apply the derived recursion operator (\ref{recur})
to the system of translations in $x$,
that is, the system of hydrodynamic type
\be
v^i_t = v^i_x,  \label{transl}
\ee
which is, obviously, Hamiltonian with the
Hamiltonian operator (\ref{op1}):
\be
v^i_t =v^i_x \equiv P^{ij}_2 {\delta H \over \delta v^j (x)},
\ \ \ H = {1 \over 2} \int \eta_{jl} v^j (x) v^l (x) dx.
\ee

Any system from the hierarchy
\be
v^i_{t_n} = \left ( R^n \right )^i_j v^j_x,
\ \ \ \ n \in {\bf Z}, \label{ierarkh}
\ee
is a multi-Hamiltonian integrable system.

In particular, any system of the form
\be
v^i_{t_1} = R^i_j v^j_x,
\ee
that is, the system of hydrodynamic
type
\bea
&&
v^i_{t_1} =
 \left ( \biggl ( \eta^{is} {\partial h^j \over \partial v^s}
+ \eta^{js} {\partial h^i \over \partial v^s} \biggr )
{d \over dx} + \eta^{is}
{\partial^2 h^j \over \partial v^s \partial v^k}
v^k_x \right ) \eta_{jl} v^l \nn\\
&&
\equiv  \left (  \eta^{is} {\partial h^j \over \partial v^s} \eta_{jk}
+  {\partial h^i \over \partial v^k}
 + \eta^{is} \eta_{jl}
{\partial^2 h^j \over \partial v^s \partial v^k} v^l \right )
v^k_x
\equiv \left ( h^i (v)+ \eta^{is} {\partial h^j \over
\partial v^s} \eta_{jl} v^l \right )_x,  \label{canon}
\eea
where $h^i(v),$ $1 \leq i \leq N,$ is an arbitrary
solution of the integrable system (\ref{(2.17)}), (\ref{(2.18)}),
is integrable.

This system of hydrodynamic type is bi-Hamiltonian
with the pair of ``canonical''
Dub\-ro\-vin--Novikov Hamiltonian operators (\ref{op1}),
(\ref{op2}):
\be v^i_{t_1} =
 \left ( \biggl ( \eta^{is} {\partial h^j \over \partial v^s}
+ \eta^{js} {\partial h^i \over \partial v^s} \biggr )
{d \over dx} + \eta^{is}
{\partial^2 h^j \over \partial v^s \partial v^k}
v^k_x \right ) {\delta H_1 \over \delta v^j (x)}, \ \ \
H_1 = {1 \over 2} \int \eta_{jl} v^j (x) v^l (x) dx, \label{eq1}
\ee
\be
v^i_{t_1} = \eta^{ij} {d \over dx}
{\delta H_2 \over \delta v^j (x)}, \ \ \
H_2 =  \int \eta_{jk} h^k (v (x)) v^j (x) dx. \label{eq2}
\ee

The next system in the hierarchy (\ref{ierarkh})
is the integrable system of hydrodynamic type
\bea
&&
v^i_{t_2} =
 \left ( \biggl ( \eta^{is} {\partial h^j \over \partial v^s}
+ \eta^{js} {\partial h^i \over \partial v^s} \biggr )
{d \over dx} + \eta^{is}
{\partial^2 h^j \over \partial v^s \partial v^k}
v^k_x \right ) \eta_{jl}
\left ( h^l (v)+ \eta^{lp} {\partial h^r \over
\partial v^p} \eta_{rq} v^q \right ) \nn\\
&&
\equiv \left ( \biggl ( \eta^{is} {\partial h^j \over \partial v^s}
+ \eta^{js} {\partial h^i \over \partial v^s} \biggr )
\left ( \eta_{jl} {\partial h^l \over \partial v^k} +
\eta_{rk} {\partial h^r \over \partial v^j} +
\eta_{rq} v^q {\partial^2 h^r \over \partial v^j
\partial v^k} \right ) \right.
\nn\\
&&
+ \left. \eta^{is}
{\partial^2 h^j \over \partial v^s \partial v^k}
\left ( \eta_{jl} h^l (v) +
\eta_{rq} v^q {\partial h^r \over \partial v^j} \right ) \right )
v^k_x.
\eea

The hierarchy of integrable systems (\ref{ierarkh})
is ``canonical'' for all bi-Hamiltonian systems
of hydrodynamic type possessing pairs of compatible
local Hamiltonian operators of hydrodynamic type.

We have also realized a completely similar explicit construction
of the corresponding class of integrable bi-Hamiltonian systems
of hydrodynamic type in the case of compatible
nonlocal Hamiltonian operators of hydrodynamic type
(see \cite{[20]}--\cite{[25]}
about the nonlocal Hamiltonian operators of hydrodynamic type).
These results will be published somewhere else.


\medskip

\begin{flushleft}
Centre for Nonlinear Studies,\\
L.D.Landau Institute for Theoretical Physics, \\
Russian Academy of Sciences\\
e-mail: mokhov@mi.ras.ru; mokhov@landau.ac.ru\\
\end{flushleft}

\end{document}